\documentclass[11pt,a4paper,leqno,twoside, headinclude]{amsart}

\title[Formality in special holonomy]{On the formality of nearly K\"ahler manifolds and of Joyce's examples in $G_2$-holonomy}

\def\auth{Manuel Amann and Iskander A.~Taimanov}
\date{December 20th, 2020}


\author{\auth}

\usepackage[english]{babel}

\usepackage{color}

\usepackage{amsmath, amssymb, amscd, amsthm}
\usepackage{stmaryrd}
\usepackage{latexsym}
\usepackage[all]{xy}
\usepackage{rotating}
\usepackage{multicol}
\usepackage{lscape}
\usepackage{wasysym}
\usepackage{longtable}
\usepackage{enumerate}

\xyoption{all}

\newtheorem{theo}{Theorem}[section]
\newtheorem{main}{Theorem}

\newtheorem*{main*}{Theorem}

\newtheorem*{mainprop*}{Proposition}

\newtheorem{mainconj}{Conjecture}

\newtheorem{defi2}[theo]{Definition}
\newtheorem*{defi2*}{Definition}
\newenvironment{defi}{\begin{defi2}\normalfont}{\end{defi2}}
\newenvironment{defi*}{\begin{defi2*}\normalfont}{\end{defi2*}}

\newenvironment{defin*}[1]{\begin{defi2*}[#1]\normalfont}{\end{defi2*}}
\newtheorem*{rem2*}{Remark}
\newenvironment{rem*}{\begin{rem2*}\normalfont}{\hfill$\boxbox$\end{rem2*}}
\newtheorem{rem2}[theo]{Remark}
\newenvironment{rem}{\begin{rem2}\normalfont}{\hfill$\boxbox$\end{rem2}}

\newtheorem*{cor*}{Corollary}

\newtheorem{conj}[theo]{Conjecture}

\newtheorem*{conj*}{Conjecture}
\newtheorem*{theo*}{Theorem}
\newtheorem*{ques*}{Question}
\newtheorem*{mi2}{Main Idea}

\newtheorem{ex2}[theo]{Example}

\newtheorem{exer2}[theo]{Exercise}

\newtheorem{alg2}[theo]{Algorithm}

\newcommand{\cc}{{\mathbb{C}}}                                     
\newcommand{\hh}{{\mathbb{H}}}                                     
\newcommand{\nn}{{\mathbb{N}}}                                     
\newcommand{\qq}{{\mathbb{Q}}}                                     
\newcommand{\rr}{{\mathbb{R}}}                                     
\newcommand{\pp}{{\mathbf{P}}}                                     
\newcommand{\s}{{\mathbb{S}}}                                      
\newcommand{\zz}{{\mathbb{Z}}}                                     
\newcommand{\SO}{{\mathbf{SO}}}                                    
\newcommand{\U}{{\mathbf{U}}}                                      
\newcommand{\SU}{{\mathbf{SU}}}                                    
\newcommand{\Sp}{{\mathbf{Sp}}}                                    
\newcommand{\G}{{\mathbf{G}}}                                      
\newcommand{\Spin}{{\mathbf{Spin}}}                                
\newcommand{\dif} {{\operatorname{d}}}                             
\newcommand{\In} {{\,\subseteq\,}}                                 
\newcommand{\Hom}{{\operatorname{Hom}}}                            
\newcommand{\id}{{\operatorname{id}}}                              
\newcommand{\APL}{{\operatorname{A_{PL}}}}                         
\newcommand{\ADR}{{\operatorname{A_{DR}}}}                         
\newcommand{\co}{\colon\thinspace}                                 

\newcommand{\comment}[1]{}                                         
\newcommand{\xto}[1]{\xrightarrow{#1}}                             
\newcommand{\hto}[1]{\overset{#1}{\hookrightarrow}}                

\newcommand{\ack}{\noindent\textbf{Acknowledgements. }}            
\newcommand{\str}{\noindent\textbf{Structure of the article. }}    

\begin{document}

\maketitle \thispagestyle{empty}


\begin{abstract}
It is a prominent conjecture (relating Riemannian geometry and algebraic topology) that all simply-connected compact manifolds of special holonomy should be formal spaces, i.e.~their rational homotopy type should be derivable from their rational cohomology algebra already---an as prominent as particular property in rational homotopy theory. Special interest now lies on exceptional holonomy $\G_2$ and $\Spin(7)$.

In this article we provide a method of how to confirm that the famous Joyce examples of holonomy $\G_2$ indeed are formal spaces; we concretely exert this computation for one example which may serve as a blueprint for the remaining Joyce examples (potentially also of holonomy $\Spin(7)$).

These considerations are preceded by another result identifying the formality of manifolds admitting special structures: we prove the formality of nearly K\"ahler manifolds.

A connection between these two results can be found in the fact that both ``special holonomy'' and ``nearly K\"ahler'' naturally generalize compact K\"ahler manifolds, whose formality is a classical and celebrated theorem by Deligne--Griffiths--Morgan--Sullivan.
\end{abstract}


\section*{Introduction}

In topology we find a beautiful synthesis of two both deeply linked yet ``competing'' theories: (co)homology and homotopy. If one is willing to forget torsion information, a similar picture is preserved, yet is it enriched by the beauty of rational homotopy theory which provides a common framework to observe the interplay of homotopy and homology ``in wildlife''. The ``rational homotopy type'' thus is elegantly encoded in so-called Sullivan models, or, respectively, in free Lie models, and contains both homotopical as well as (co)homological invariants. This immediately prompts one of the most considered meta-questions in the field: when is the homotopy type no more complicated than the cohomology type, when is all homotopical information already contained in the cohomology algebra, i.e.~when is the space \emph{formal}?

It is one of the mysteries of the manifold interlinks between topology and (Riemannian) geometry that this abstract concept of formality suddenly seeks attention in the form of analytically defined properties of Riemannian structures. More concretely, formality holds for symmetric spaces (and their generalizations), it is satisfied by K\"ahler manifolds, and it is conjectured (and partly confirmed) to hold true for Riemannian manifolds of special holonomy (see Sections \ref{sec01}, \ref{sec04}).

With this article we want to cast more light on the depicted ``opaque interweavings'' of geometry and topology with our two main theorems.
\begin{main}\label{theoA}
Closed simply-connected nearly K\"ahler manifolds are formal.
\end{main}
Nearly K\"ahler manifolds are one way to generalize K\"ahler manifolds and we point the reader to Section \ref{sec02} for their definition. In order to prove this result, it suffices to bring together known classification results for nearly K\"ahler manifolds on the one hand, and, on the other hand, formality results for the respective ``building blocks'' in these classifications.

\bigskip

It is a prominent conjecture that all simply-connected manifolds of special holonomy are formal. With the formality of K\"ahler manifolds (including hyperk\"ahler ones) and of positive quaternion K\"ahler manifolds established the attention focuses on $\G_2$- and $\Spin(7)$-manifolds. In both cases there are several simply-connected closed examples, the first and most prominent ones established by Joyce. Until recently (see \cite{Tai18}) not even their cohomology rings were known. We extend the analysis of the geometric construction of the example in this article via means from differential topology and rational homotopy. More precisely, we use the intersection homology of concrete submanifolds in order to build a Sullivan model over the reals up to a certain degree. For a typical example of a Joyce manifold we consider the simplest one, which we denote simply by $M$
(see \cite[p. 349, Example 1]{Joyce-g2-2}). Note that this example, however, already carries the main features of the construction principles and can serve as a model case (see Remark \ref{rem02}). Consequently, the concept of $s$-formality allows us to derive

\begin{main}\label{theoB}
The Joyce example $M$ 
of exceptional holonomy $\G_2$ is formal.
\end{main}

The manifold in question is described in detail in \cite{Tai18} and in Section \ref{sec03}.

This approach and reasoning can be considered a blueprint result in order to analyze and establish the formality of the remaining Joyce examples (in $\Spin(7)$ holonomy as well) (cf.~Remark \ref{rem02}).
The argumets are purely topological and are not related directly to the existence of a metric with special holonomy. However we would like to stress that following prominent conjecture (see for example \cite[Conjecture 2, p.~2049]{AK} guides our work with respect to special holonomy.
\begin{conj}
A simply-connected compact Riemannian manifold of special holonomy is a
formal space.
\end{conj}

\bigskip

\str In Section \ref{sec01} we provide necessary concepts from rational homotopy theory. In particular, we provide a certain background on formality. Section \ref{sec02} is devoted to compiling necessary results and to the proof of Theorem \ref{theoA}. In Section \ref{sec04} we comment on the importance of formality within special holonomy, we review the (generalized) Kummer construction underlying the examples by Joyce, and we recall the concepts from intersection homology which we use on the side of differential topology. Next, in Section \ref{sec03} we provide the proof of Theorem \ref{theoB}. This illustrates the general method of deriving formality of such examples.

\ack
The authors thank Johannes Nordstr\"om for recalling the Bianchi--Massey tensor (see Remark \ref{rem03}) as well as an anonymous referee for helping to improve the presentation of the article. We thank Luc\'ia Mart\'in-Merch\'an for pointing out a wrongly chosen differential form.
We thank the referee for helpful comments.

The first named author was supported both by a Heisenberg grant and his research grant AM 342/4-1 of the German Research Foundation; he was moreover a member of the DFG Priority Programme 2026 and thanks the University of Augsburg for their continued support.
The second named author was supported by RFBR (grant 18-01-00411-a) and
was performed according to the Government research assignment for IM SB RAS, project FWNF-2022-0004.

\section{Formality in Rational Homotopy Theory}\label{sec01}

We cannot and do not intend to provide an introduction to rational homotopy theory here. For this we point the reader to the textbooks \cite{FHT01}, and \cite{FOT08}.
Let us just briefly recall some concepts necessary for this article.

A most important tool in the theory are so-called \emph{(minimal) Sullivan models} of commutative differential graded algebras, respectively of nilpotent or simply-connected topological spaces $X$, which are encoded via the algebras of polynomial differential forms $(\APL(X),\dif)$ (see \cite[Chapters 10 and 12]{FHT01}). Indeed, the rational homotopy type of any such space is captured via a quasi-isomorphism (i.e.~a morphism of differential graded algebras inducing an isomorphism on cohomology)
\begin{align*}
(\Lambda V,\dif) \xto{\simeq} (\APL(X),\dif)
\end{align*}
from a minimal Sullivan algebra $(\Lambda V,\dif)$. For path-connected spaces such a model can always be constructed and the algebra $(\Lambda V,\dif)$ is unique up to isomorphism. For example, the fact that it encodes the rational homotopy type of $X$ is reflected by the fact that up to duality the underlying graded vector space $V$ over $\qq$ is gradedly isomorphic to $\Hom(\pi_*(X),\qq)$, the rational homotopy groups of $X$, provided that $X$ is a nilpotent space.

It is one of the central questions in rational homotopy theory to understand the nature, behaviour and examples of \emph{formality}.
\begin{defi}
A path-connected topological space $X$ is called \emph{formal} if already its rational cohomology algebra $H^*(X;\qq)$ (endowed with the trivial differential) encodes its rational homotopy type, i.e.~if there is a chain of quasi-isomorphisms of cochain algebras over the rationals
\begin{align*}
(\APL,\dif)\xleftarrow{\simeq} (A_1,\dif) \xto{\simeq} \dots \xleftarrow{\simeq} (A_k,\dif) \xto{\simeq} (H^*(X),0)
\end{align*}
\end{defi}

So far, the algebras we consider were over the rational numbers. Clearly, one obtains the analog results passing to algebras over the reals, thus dealing with the \emph{real homotopy type}. On smooth manifolds one may observe that the commutative differential graded algebra of smooth differential forms $(\ADR(M),\dif)$ together with the exterior derivative encodes the real homotopy type. Usually, several rational homotopy types may fall into one real homotopy type. However, from \cite[Corollary 6.9, p.~265]{HS79} we recall that a space is formal over $\qq$ if and only if it is formal over any field extension, i.e.~in particular over $\rr$ (with the analog definition replacing polynomial differential forms by smooth differential forms). This result will allow us to draw on computations with smooth differential forms in order to establish formality.

\bigskip

There are several prominent examples of formal spaces. Since products of harmonic forms on compact symmetric spaces are harmonic again, these spaces are ``geometrically formal'' and formal, in particular, due to the Hodge decomposition. We remark that this is one of the few cases (cf.~Theorem \ref{theoDGMS}, see Remark \ref{rem01}) in which formality can be derived by an analytic property.

There is the larger class of so-called \emph{generalized symmetric spaces}. These are Riemannian manifolds equipped with a \emph{(regular) s-structure}, i.e.~at each point $x\in M$ there exists an isometry $s_x\co M\to M$ with isolated fixed point $x$ and satisfying that
\begin{align*}
s_xs_y=s_{s_x(y)}s_x
\end{align*}
In case there exists a $k\in \nn$ such that for all $x\in M$ it holds that $s_x^k=\id$ and $s_x^l\neq \id$ for $l<k$, then $k$ is the \emph{order} of the s-structure. The smallest order $k$ of any s-structure on $(M,g)$ is the order of the generalized symmetric space. We then  call $(M,g)$ a \emph{$k$-symmetric space}. We remark that the given definition is not the most general one used in the literature, but suffices for our purposes.

From the main results in \cite{S}, \cite{KT}, \cite{GN16} we cite that $k$-symmetric spaces are formal (but not necessarily geometrically formal).

\bigskip

This article is particularly motivated by another highly important class of formal spaces, namely compact K\"ahler manifolds.
\begin{theo}[\cite{DGMS}]\label{theoDGMS}
Compact K\"ahler manifolds are formal.
\end{theo}

Moreover, since, due to the K\"unneth formula, the cohomology of a finite product space is the tensor product of the respective cohomology algebras, and since a Sullivan model of such a tensor product is the tensor product of the respective models, it is immediate that a product of formal spaces is again formal---actually, it is formal if and only if so are both spaces.

One-point unions of formal topological spaces and connected sums of formal manifolds are formal again.


\section{Formality of nearly K\"ahler manifolds}\label{sec02}

There are two natural generalizations of K\"ahler manifolds: almost K\"ahler manifolds
and nearly K\"ahler manifolds. These notions are orthogonal in the sense that if a manifold is
almost K\"ahler and nearly K\"ahler at the same time, then it is K\"ahler (see \cite{GH}).

We recall that a manifold is called {\it almost Hermitian} if it is endowed with
an almost complex structure $J$ and Hermitian metrics in the fibers of the tangent bundle.
It is assumed that $J$ and the metric are smooth.

\bigskip

An {\it almost K\"ahler} manifold is an almost Hermitian manifold for which the smooth $2$-form on the fibers of the tangent bundles
\begin{align}
\label{symp}
\omega(X,Y) = \langle JX,Y \rangle
\end{align}
is closed. Here $\langle \cdot, \cdot \rangle$ is the Hermitian metric.

It is easy to notice that these are exactly symplectic manifolds for which $\omega$ is a symplectic form.
There was a conjecture that simply-connected compact symplectic manifolds are formal. For symplectic nilmanifolds it was known to be wrong: for instance, the Kodaira--Thurston manifold $M_{KT} = M_H \times S^1$, where $M_H$ is a nilpotent Heisenberg $3$-manifold, is not formal due to a nontrivial triple Massey product in $M_H$.   In \cite {BaT98,BaT00} the formality conjecture for simply-connected symplectic manifolds of dimension $ \geq 10$ was disproved as follows. $M_{KT}$ is symplectically embedded
into $\cc \pp^n, n \geq 5$, and after the symplectic blow up of $\cc \pp^n$ along $M_{KT}$ the nontrivial Massey product in $M_{KT}$ gives rise to a nontrivial Massey product in the resulting manifold. For eight-dimensional manifolds this conjecture was disproved later by other means (see \cite{FM08}).

\bigskip

A {\it nearly K\"ahler manifold} is an almost Hermitian manifold such that
the $(2,1)$-tensor $\nabla J$ is skew symmetric:
\begin{align}
\label{nk}
(\nabla_X J)(X) = 0 \ \ \ \mbox{for all $X \in TM$}.
\end{align}

This notion was introduced in \cite{Gray70} where it was mentioned that there are plenty of examples of nearly K\"ahler manifolds which are not K\"ahler. Among these are the $6$-sphere with the canonical almost complex structure and the round metric, and the quotients $G/K$ where $G$ is a compact semisimple Lie group, and $K$ is the fixed point set of an automorphism $\sigma: G \to G$ of order three: $\sigma^3=1$. These are $3$-symmetric spaces in the terminology of Section \ref{sec01}.

We remark that it is a prominent problem if the $6$-sphere admits a complex structure; clearly, it is certainly not K\"ahler.

Contrast the definition of nearly K\"ahler manifolds with the definition of a K\"ahler manifold wherein, in particular, it is required that $J$ be a genuine complex structure, i.e.~$\nabla J=0$ and $J$ is parallel. Thus, clearly, nearly K\"ahler manifolds are obvious generalizations.

It appears that, in contrast to \eqref{symp}, Condition \eqref{nk} is very rigid.

\bigskip

Let us recall that a nearly K\"ahler manifold is called {\it strict} if
$$
\nabla_X J \neq 0 \ \ \ \mbox{for $X \neq 0$.}
$$
It was proved in \cite{Kirichenko} (see also \cite{Nagy02b}) that  a complete and simply-connected
nearly K\"ahler manifold $M$ is represented as a Riemannian product
\begin{align}\label{eqnnKsplit}
M = M_1 \times M_2
\end{align}
of a K\"ahler manifold $M_1$ and a strict nearly K\"ahler manifold $M_2$ (one of the factors may be a point).

The classification of complete simply-connected strict nearly  K\"ahler manifolds is given by Nagy's Theorem in \cite{Nagy02a}: every such manifold is a Riemannian product whose factors belong to one of three classes, namely

\begin{enumerate}
\item
$6$-dimensional nearly K\"ahler manifolds,

\item
homogeneous nearly K\"ahler manifolds of certain types,

\item
twistor spaces over quaternionic K\"ahler manifolds with positive scalar curvature,
endowed with the canonical nearly K\"ahler metric.
\end{enumerate}

Moreover, Nagy derived from this classification the proof of the Wolf--Gray conjecture for manifolds of dimension $\geq 8$. This conjecture states that every  homogeneous nearly K\"ahler manifold is a $3$-symmetric space with the canonical almost complex structure, and its proof was completed by Butruille who confirmed it for $6$-dimensional manifolds in \cite{But}.

\bigskip

In order to prove the formality of nearly K\"ahler manifolds, we recall the following observations from Section \ref{sec01}. First of all, from Theorem \eqref{theoDGMS} we recall the formality of compact K\"ahler manifolds. Next, a product is formal if (and only if) so are its factors. Hence, in view of Decomposition \eqref{eqnnKsplit} it remains to prove the formality of strict nearly K\"ahler manifolds, i.e.~we have to discuss formality for the three possible types of factors given above.

\begin{enumerate}
\item all simply-connected closed $6$-dimensional manifolds are formal for purely algebraic reasons (see \cite{NM77});

\item $3$-symmetric manifolds and moreover all $k$-symmetric spaces $G/H$ with $G$, $H$ compact connected are formal by \cite{GN16,KT,S};

\item positive quaternion K\"ahler manifolds and twistor spaces over them are formal (see \cite{DGMS}, \cite{AK}).
\end{enumerate}

Let us quickly reflect upon and justify the third item. Quaternion K\"ahler manifolds are Einstein. If their scalar curvature is positive, they are called \emph{positive quaternion K\"ahler manifolds}. In particular, they are simply-connected then (see \cite[Theorem 6.6, p.~163]{Sal82}).

Quaternion K\"ahler manifolds admit so-called \emph{twistor fibrations} $\cc\pp^1\hto{}E\to M$ with complex contact Einstein \emph{twistor spaces} $E$. The complex structure is constructed using the quaternionic structure of the base and the complex structure of $\cc\pp^1$---see \cite[Theorem 4.1, p.~152]{Sal82}. (The simplest example of this is given by $\cc\pp^1\hto{}\cc\pp^{2n+1}\to \hh\pp^n$.)  Due to \cite[Theorem 6.1, p.~158]{Sal82} twistor spaces $E$ over positive quaternion K\"ahler manifolds actually admit a K\"ahler structure. Hence they are formal. (Note that in \cite{AK} this formality was used to derive the formality of positive quaternion K\"ahler manifolds as well.) We should remark that this K\"ahler structure on $E$ is a priori independent of the nearly K\"ahler structure, which justifies that, endowed with the nearly K\"ahler structures, the manifolds $E$ are actually strict nearly K\"ahler manifolds.

Since a product of formal manifolds is formal, this then yields

\begin{theo}
Simply-connected closed nearly K\"ahler manifolds are formal.
\end{theo}

\begin{rem}\label{rem01}
The formality of K\"ahler manifolds is derived from the so-called $dd^c$-lemma.
It would be interesting to derive the formality of nearly K\"ahler manifolds straightforwardly from the analytical condition \eqref{nk} without checking one by one special cases as done above.
\end{rem}


\section{Formality and special holonomy}\label{sec04}

\subsection{The formality problem for manifolds with special holonomy}
There is another way to generalize Theorem \ref{theoDGMS} by Deligne--Griffiths--Morgan--Sullivan on the formality of K\"ahler manifolds. This is based on the observation that K\"ahler manifolds are exactly those Riemannian manifolds which have $\U(n)$-holonomy.

By Berger's Theorem, the holonomy group of an $n$-dimensional oriented manifold which is neither locally reducible (i.e.~not locally a product) nor locally symmetric
either coincides with the entire group $\SO(n)$ or has \emph{special holonomy} from the following list:
$\U\left(\frac{n}{2}\right)$ (K\"ahler),
$\SU\left(\frac{n}{4}\right)$ (Calabi--Yau),
$\Sp\left(\frac{n}{4}\right)$ (hyperk\"ahler),
$\Sp\left(\frac{n}{4}\right)\Sp(1)$ (quaternion K\"ahler),
$\G_2$ ($n=7$), and
$\Spin(7)$ ($n=8$).

Since $\Sp\left(\frac{n}{4}\right) \subset \U\left(\frac{n}{2}\right)$ and
$\SU\left(\frac{n}{2}\right) \subset \U\left(\frac{n}{2}\right)$, manifolds whose holonomy belongs to the
first three special types are K\"ahler. By Theorem \ref{theoDGMS}, they are formal.

As we recalled in Section \ref{sec02} quaternion K\"ahler manifolds are Einstein. If they are simply-connected and have vanishing scalar curvature, they are hyperk\"ahler, K\"ahler and formal, in particular. Positive quaternion K\"ahler manifolds (those with positive scalar curvature) are simply-connected due to \cite[Theorem 6.6, p.~163]{Sal82} and formal due to \cite{AK}. If they have negative Einstein constant, they are called \emph{negative quaternion
K\"ahler manifolds} (of which we do not know many interesting examples).

Hence we are left with three open cases:

\begin{itemize}
\item
negative quaternion K\"ahler manifolds,

\item $\G_2$-manifolds,

\item $\Spin(7)$-manifolds.
\end{itemize}

Let us recall that the first examples of closed manifolds with holonomy $\G_2$ and $\Spin(7)$ were found in the 1990s by Joyce (see \cite{Joyce-g2-1,Joyce-g2-2,Joyce-spin7}). He used for that the generalized Kummer construction. Since then other methods of constructing such examples were introduced
(see \cite{Joyce99,Kovalev,KovalevLee,CHNP,JK, Nor18}).

We would like to consider the formality problem for manifolds constructed with this method by Joyce.

\subsection{The Kummer construction and its generalizations}

The original Kummer method deals with the involution $\sigma\co T^4 \to T^4$ which in the linear coordinates
$x = (x_1,x_2,x_3,x_4)$ on the four-torus $T^4 = \rr^4/\zz^4$ takes the form $\sigma(x)=-x$.
We recall that these coordinates are defined modulo integers and that the involution has exactly $16$ fixed points,
i.e.~the half-periods, which are defined by the condition $2x \in \zz^4$. The singular set  of $T^4/\langle \sigma \rangle$ consists of $16$ conic points whose neighborhoods have the form $B^4/\{\pm 1\}$ where $B^4$ is the open $4$-dimensional disc. Each singularity is resolved via the mapping $T^\ast \cc \pp^1 \to \cc^2/\{\pm 1\} \supset B^4/\{\pm 1\}$, i.e., every copy of $B^4/\{\pm 1\}$ is replaced by an open neigborhood $U$ of the zero section of the fiber bundle $T^\ast \cc\pp^1 \to \cc\pp^1$.
The space $T^\ast\cc\pp^1$ admits the Ricci-flat asymptotically locally Euclidean (ALE) metric called the Eguchi--Hanson metric. The open manifold $T^4/\langle\sigma\rangle$ is flat. By perturbations these metrics are made compatible and form a Ricci-flat metric on the resulting four-manifold, which is homeomorphic to a $K3$ surface. Since the preimage of zero under the mapping $U \to B^4/\{\pm 1\}$ is a two-sphere,
the resolution of each singularity increases $b_2(T^4/\langle \sigma \rangle)$, the second Betti number of  $T^4/\langle \sigma \rangle$, by one.

\bigskip

In  \cite{Joyce-g2-1,Joyce-g2-2,Joyce-spin7} this construction is generalized as follows. Joyce considers
an action of a finite group $\Gamma$ on a flat torus $T^7$ or $T^8$ such that the singular set $T/\Gamma$ consists of finitely many orbifolds. In fact, in these articles finitely many examples for which every orbifold singularity is resolved by an analogue of the Kummer construction are considered.

The simplest case is the orbifold $T^3$ with the neighborhood $T^3 \times B^4/\{\pm 1\}$. Such a singularity is  resolved by applying the Kummer construction to every fiber of the fibration $T^3 \times B^4/\{\pm 1\} \to T^3$.   Therewith, $T^3 \times B^4/\{\pm 1\}$ is replaced by $T^3 \times U$ where $U$ is a neighborhood of the zero section of the bundle $T^\ast \cc\pp^1 \to \cc\pp^1$, the resolution map $\pi\co U \to B^4/\{\pm 1\}$
is  a diffeomorphism on  the preimage of the punctured disc $(B^4 \setminus \{0\})/\{\pm 1\}$, and
$\pi^{-1}(0)$ is the zero section of the bundle   $T^\ast \cc\pp^1 \to \cc\pp^1$. Therefore,
the preimage of $T^3 \times \{0\}$ under the resolution map is a product $T^3 \times \cc\pp^1$, and we can easily describe the contribution of the resolution to the Betti numbers and the intersections of the appearing cycles.

For seven-dimensional manifolds the other possible singularities are
\begin{itemize}
\item
either $T^3/\langle \tau\rangle$, where in linear coordinates the involution $\tau$ is given by the formula
$\tau(x_1,x_2,x_3) = (x_1+\frac{1}{2}, -x_2,-x_3)$, and the neighborhoods look like $(T^3 \times B^4/\{\pm 1\})/\zz_2$,

\item
or $(T^3 \times B^4/\zz_3)/\zz_2$, and $\s^1$ for which there are three possible types of neighborhoods: $\s^1 \times B^2/\zz_3, (\s^1 \times  B^2/\zz_3)/\zz_2$, and $\s^1 \times B^6/\zz_7$. In the last three cases the Eguchi--Hanson metric is replaced by other explicitly described Ricci-flat ALE metrics.
\end{itemize}

For eight-dimensional manifolds five model types of singularities are explicitly described in \cite[Section 3]{Joyce-spin7}. In the seven-dimensional case all model singular orbifolds are pairwise-nonintersecting. In the eight-dimensional case one of the model singularities appears as the intersection of singular orbifolds of two other types.

\subsection{The intersection homology ring}
As a general reference for this section we point the reader to \cite[Section 6, p.~53]{BT82}.

The homology groups $H_\ast(M;\rr)$ of manifolds obtained by the generalized Kummer construction
are generated by cycles from $H_\ast(T^{7\, (\mathrm{or}\, 8)};\rr)$ invariant under $\Gamma$ and by the explicitly described cycles induced by the resolutions.

Let us assume that the manifold $M^n$ is oriented. It is known that if two homology cycles $u \in H_k(M^n;\rr)$ and $v \in H_l(M^n;\rr)$ are realized by oriented submanifolds $X^k$ and $Y^l$ and intersect transversally, then their intersection $X \cap Y$ realizes the cycle $w$ such that
$$
Dw = Du \cup Dv
$$
where
$$
D: H_\ast(M^n;\rr) \to H^\ast(M^n;\rr)
$$
is Poincare duality (see, for instance, \cite[Section 6]{BT82}). Here the orientation of $X \cap Y$ is defined as follows.
Let $(e_1,\dots, e_{k+l-n})$ be a basis in the tangent space to $X\cap Y$ at some point $x$, and let
$(e_1,\dots, e_{k+l-n}, e^\prime_1,\dots,e^\prime_{n-l})$ and $(e_1,\dots , e_{k+l-n}, e^{\prime\prime}_1,\dots,e^{\prime\prime}_{n-k})$ be positively oriented bases in the tangent spaces to $X$ and $Y$ at the same point. Now we are left to define the orientation of $X \cap Y$.  If the basis $(e_1,\dots, e_{k+l-n},e^\prime_1,\dots,e^\prime_{n-l},e^{\prime\prime}_1,\dots,e^{\prime\prime}_{n-k})$ at the tangent space to $M^n$
is positively oriented, we assume that $(e_1,\dots, e_{k+l-n})$ is positively oriented in the tangent space to $X \cap Y$. Otherwise, we assume that the orientation of $X \cap Y$ is determined by the basis
$(-e_1,e_2,\dots, e_{k+l-n})$.

Moreover, for every such cycle $u$ its Poincare dual can be realized by a closed form $\omega$ supported on an open neighborhood of the corresponding submanifold $X$ (see \cite{BT82}). Such a neighborhood can be taken arbitrarily small.

This way we obtain the homology intersection ring which is dual to the cohomology ring.
It was originally defined by Lefschetz for algebraic varieties and later his sketch of the definitions
for all oriented manifolds was rigorously realized in \cite{GP}. However, if all generators are realized by submanifolds this ring is quite clear.

From the generalized Kummer construction it is possible to write down explicitly all (additive) generators of the homology, to realize them by submanifolds which are pairwise transversally intersecting, and to
explicitly describe the ring structure. For the simplest example, the Joyce manifold with all singularities of the form $T^3 \times B^4/\{\pm 1\}$ this was done in \cite{Tai18}. In fact, in \cite{Tai18} there was presented an algorithm for computing case by case the (real) cohomology rings for all of Joyce's manifolds. For another example, this was realized in the diploma work of I.V.~Fedorov, a student of
the second named author (I.A.T.) (see \cite{Fedorov}).

The example from \cite{Tai18} is used in the next section. The reader should understand and interpret this as a blueprint of how to prove the formality of Joyce's
$\G_2$-manifolds. We suggest to do this following the same method.

\begin{rem}\label{rem02}
We think that after some technical modifications---keeping also in mind that the orbifold singularities
in this case may intersect---the same method can be applied to prove the formality of not only the Joyce examples of holonomy $\G_2$, but also of Joyce's $\Spin(7)$-manifolds.
\end{rem}


\section{Formality of the Joyce examples}\label{sec03}

\subsection{Strategy and Preliminaries}
Denote by $M$ the Joyce example from \cite{Tai18} with holonomy $\G_2$. We want to prove that $M$ is formal in the sense of Rational Homotopy Theory. Since formality does not depend on the field extension of $\qq$ (see \cite[Corollary 6.9, p.~265]{HS79}), we work with real coefficients, i.e.~with usual deRham differential forms $\ADR(M)$ and with minimal Sullivan models and cohomology over the reals.

We recall \cite[Theorem 3.1, p.~157]{FM05}, namely
\begin{theo}\label{theo01}
Let $M$ be a connected and orientable compact differentiable manifold
of dimension $2n$ or $2n-1$. Then $M$ is formal if and only if it is $(n-1)$-formal.
\end{theo}
The theorem draws on the concept of $s$-formality originally introduced in \cite[p.~150]{FM05}. For the convenience of the reader we directly specify this to $s=3$ below.

Hence, in order to prove the formality of $M$ it suffices to show that $M$ posseses some $3$-formal minimal model $(\Lambda V,\dif)$. That is we shall construct a minimal model $(\Lambda V,\dif)$ of $M$ with the property that
\begin{itemize}
\item
there is a homogeneous splitting $V=C\oplus N$ with $\dif C=0$,
\item
$\dif\co N\to \Lambda V$ is injective,
\item
and any closed form in the ideal $I(N^{\leq 3})=N^{\leq 3} \cdot \Lambda V^{\leq 3}$ is exact in $(\Lambda V,\dif)$.
\end{itemize}

This minimal model will be derived from a special model $(\mathcal{A},\dif)$ of $M$ which we shall construct in the following.

From \cite[Theorem 1, p.~9]{Tai18} we cite the structure of the rational cohomology ring. It is the quotient of
\begin{align*}
\qq[c_{\delta i}, c_{\delta i j}', t_\delta', t_i']_{1\leq i\leq 4, 1\leq j\leq 3, 1\leq \delta \leq 3} \otimes  \Lambda\langle t_\delta, t_i, c_{\delta i j}, c_{\delta i}'\rangle_{1\leq i\leq 4, 1\leq j\leq 3, 1\leq \delta \leq 3}
\end{align*}
with degrees
\begin{align*}
&\deg c_{\delta i}=2, \deg c_{\delta i j}=\deg t_\delta=\deg t_i=3, \\&\deg c_{\delta i j}'=\deg t_\delta'=\deg t_i'=4, \deg c_{\delta i}'=5
\end{align*}
by the relations
\begin{align*}
c_{\delta i}c_{\delta i}'=-2[M], \ c_{\delta i j}c_{\delta i j}'=-2[M], \ t_\delta t_\delta'=8[M], \ t_it_i'=8[M], \
c_{\delta i}^2=-2t_\delta',
\end{align*}
with $[M]$ the fundamental class and with all other products of generators vanishing.

\subsection{Construction of the minimal model}
As was done in \cite{Tai18}, up to duality, we now choose representing embedded submanifolds of complementary dimensions. The manifold $M$ is obtained from $T^7$ by resolving singularities of $T^7/\Gamma$, where the group $\Gamma = \zz^3_2$ is generated by three involutions. The singularity set in $T^7/\Gamma$ consists of $12$ three-dimensional tori $T_{\delta i}, 1\leq \delta \leq 3, 1 \leq i \leq 4$,
with neighborhoods of the form $T^3 \times B^4/\{\pm 1\}$. These tori split into three families corresponding to fixed point sets of three generating involutions numerated by $\delta$. The resolution gives rise to the submanifolds $(C_{\delta i}')^5=T_{\delta i}\times C_{\delta i}$ with $[(C_{\delta i}')]\in H_{5}(M)$ which are dual to $c_{\delta i}\in H^2(M)$ and embedded into $M$ via $\iota_{\delta i} \co C'_{\delta i}\hto{} M$. The $C_{\delta i}$ are pairwise disjoint $\cc\pp^1$s. The products of one- and two-dimensional cycles in $T_{\delta i}$ with two-spheres $C_{\delta i}$ realize the cycles $C_{\delta ij}$ and $C^\prime_{\delta ij}$. There are also
three-dimensional cycles $T_\delta$ and $T_i$ realized by tori invariant under $\Gamma$ and the
four-dimensional cycles  $T^\prime_\delta$ and $T^\prime_i$ induced by tori complemented to
$T_\delta$ and $T_i$. For instance, $T^\prime_\delta$ is diffeomorphic to a $K3$ surface.
The tori $T_{\delta i}$ are homologous to $T_\delta$ for all $\delta$ and $i$.

In this notation the cycles $[Z]$ and $[Z^\prime]$ are Poincar\'e dual to the cocycles $z^\prime $ and $z$.

From \cite[Proposition 6.24, p.~67]{BT82} we recall that the homology class represented by a closed submanifold $\iota\co  S^k\hto{} M^n$ in cohomology corresponds to the Thom class of its normal bundle, which we usually identify with a tubular neighborhood. Since the Thom class is compactly supported we may extend it to all of $M$.

Moreover, this Thom class identifies with the uniquely determined class $[\nu_S]\in H^{n-k}(M)$ satisfying
\begin{align}\label{eqn01}
\int_S \iota^* \omega=\int_M\omega\wedge \nu_S
\end{align}
for all $\omega$ with $[\omega]\in H^k(M)$ (see \cite[p.~67, (5.13), p.~51, (5.21), p.~65]{BT82}). Transversal intersection of the representing submanifolds corresponds to the cup product of their duals in cohomology (see \cite[p.~69]{BT82}). By the ``localization principle'' (see \cite[Proposition 6.25, p.~67]{BT82}, without restriction, 
a cohomology class dual to $[S]$ can be represented by a form with support contained in an arbitrary tubular neighborhood of $S$. This will play an important role in the arguments to come.

\bigskip

We construct a ``small'' sub-algebra of the deRham forms $\ADR(M)$ which is weakly equivalent to the latter.
We thus define the following commutative connected cochain algebra $(\mathcal{A},\dif)$ over $\rr$. We take $(\mathcal{A},\dif)$ to be the commutative differential graded subalgebra of $(\mathcal{A},\dif) \hto{\iota}\ADR(M)$ (together with the restriction of the usual exterior differential $\dif$) generated by the following graded subspaces $\mathcal{A}^i$ of smooth differential forms on $M$.
\begin{align*}
\mathcal{A}^0&:=\rr \\
\mathcal{A}^1&:=0 \\
\mathcal{A}^2&:=\langle c_{\delta i} \rangle_{\delta, i}\\
\mathcal{A}^{3}&:=\langle n_{\delta k}\rangle_{\delta, 2\leq k\leq 4} \oplus \langle n_\delta\rangle \oplus \langle t_\delta, t_i, c_{\delta i j}\rangle_{\delta, i, j}\\
\mathcal{A}^{\geq 4}&:=\ADR^{\geq 4}.
\end{align*}

It sounds reasonable to take for $c_{\delta i}$ representatives of the Thom classes which are supported in certain tubular neighborhoods $\nu_{\delta i}$ of the $C_{\delta i}'$ (hence representing the respective cohomology classes in $H^*(M)$ we cited above). Therewith one has to choose these tubular neighborhoods $\nu_{\delta i}$ of the $C'_{\delta i}$ small enough to guarantee that all pairwise intersections $\nu_{\delta_1 i_1}\cap \nu_{\delta_2 i_2}=\emptyset$ are disjoint whenever $\delta_1\neq \delta_2$ (which is possible, since the $C_{\delta i}$ are pairwise disjoint). However, in order to prove that $7$-dimensional Massey products vanish we construct the following representatives of these elements.


\bigskip

{\noindent \sc The construction of $c_{\delta i}$ and $n_{\delta k}$.}

Let us construct the elements $c_{\delta i}$ and $n_{\delta k}$.
We recall that the manifold $M$ is obtained by blow-up of singularities
of type $T^3 \times (\cc^2/\pm 1)$ of the quotient space $T^7/\Gamma$.

The singularities correspond to $48$ tori which are fixed by different involutions from $\Gamma$. These tori are disjoint and split into three families which correspond to the involutions $\alpha, \beta$, and
$\gamma$ that generate $\Gamma$. On every such  family the group $\Gamma$ acts by translations, every orbit of $\Gamma$ consists of $4$ tori and different orbits are mapped into each other by translations by half-periods of $T^7 = \rr^7/\zz^7$.

Since the group $\Gamma$ consists of involutions of the form
\begin{align}
\label{comm}
(x_1,x_2,\dots,x_7) \to (\pm x_1 + c_1, \pm x_2 + c_2, \dots, \pm x_7 + c_7)
\end{align}
with constants $c_1,\dots,c_7$ all translations by half-periods commute with elements of $\Gamma$.

Every such torus $T^3$ gives rise to a singularity in $T^7/\Gamma$ due to a certain element $g \in \Gamma$ which, in convenient linear coordinates $(y_1,\dots,y_n)$ on $T^3$, acts near $T^3$
as follows
$$
(y_1,\dots,y_3,y_4,\dots,y_7) \to (y_1,\dots,y_3,-y_4,\dots,-y_7),
$$
where $T^3$ is distinguished by the equation $y_4=\dots=y_7=0$.
To every such torus $T^3$ there exist neighborhoods $T^3 \subset T^3 \times U_0 \subset T^3 \times U \subset T^7 \subset T^3 \times U_1$ such that
\begin{itemize}
\item
for different tori these neighborhoods are disjoint;
\item
if tori from two orbits of $\Gamma$ are mapped into each other by a translation by a half-period then the corresponding neighborhoods are also mapped into each other by the same translation;
\item
in the local coordinates $\{y_k\}$ the sets $U_0, U$, and $U_1$ are given
by neighborhoods of the origin in $\rr^4$ with the coordinates $y_4,\dots,y_7$ and we assume that $U_0, U$ and $U_1$ are invariant under the reflection $(y_4,\dots,y_7) \to (-y_4,\dots,-y_7)$.
\end{itemize}

Now let us construct the generators $c_{\delta i}$ and $n_{\delta k}$ of the minimal model.

We recall the explciit form of the involutions $\alpha, \beta$, and $\gamma$:
$$
\alpha: (x_1, x_2, x_3,x_4,x_5,x_6,x_7) \to
(x_1,x_2,x_3,-x_4,-x_5,-x_6,-x_7),
$$
$$
\beta: (x_1, x_2, x_3,x_4,x_5,x_6,x_7) \to
(-x_1,- x_2, x_3,-x_4,\tfrac{1}{2}-x_5,x_6,x_7),
$$
$$
\gamma:
(x_1, x_2, x_3,x_4,x_5,x_6,x_7) \to
(-x_1,x_2,-x_3,\tfrac{1}{2}-x_4,x_5,\tfrac{1}{2}-x_6,x_7).
$$
Without loss of generality we assume that $\delta=\alpha$.
The involution $\alpha$ acts on the torus $T^7 = \rr^7/\zz^7$ and its fixed point set consists of
$16$ tori of the form
$(x_1,x_2,x_3,a_4,a_5,a_6,a_7)$ with $a_k \in \{0,1/2\}$, $k=4,\dots,7$.
The group $\Gamma/\zz_2\langle\alpha\rangle = \zz_2 \oplus \zz_2$ acts on these tori effectively and splits them into four orbits. Every such orbit induces a singularity in $T^7/\Gamma$, and another four singularities similarly correspond to the involutions $\beta$ and $\gamma$ each. By construction, \cite{Joyce-g2-2},
(see the more detailed description in \cite{Tai18}) the orbits of different involutions are
disjoint.

The orbifold $T^7/\langle \alpha \rangle$ is a product of a three-torus with the coordinates $x_1,x_2,x_3$ and the Kummer surface $T^4/\langle \sigma \rangle$ where $\sigma: (x_4,x_5,x_6,x_7) \to (-x_4,-x_5,-x_6,-x_7)$.
The Kummer surface containts $16$ neighborhoods, of the singular tori, of the form
$$
T^3 \times U/\langle \sigma \rangle
$$
and, since $U/\langle \sigma \rangle \subset \cc^2/\{\pm 1\}$, we take the blowups of these domains which desingularize the Kummer surface into a K3 surface $K_3$.

Therefore, after such a desingularization we
deform $T^7/\langle \alpha\rangle$ into the product
$$
\widetilde{T^7/\langle \alpha \rangle} = \widetilde{T^3 \times U/\langle \sigma \rangle} = T^3 \times K_3.
$$

To obtain $M$ we have to act on $T^3 \times K_3$ by the group $\langle \beta,\gamma\rangle = \zz_2 \oplus \zz_2$ generated by $\beta$ and $\gamma$ and desingularize the singularities. Since $\beta$ and $\gamma$ act on the fixed point tori of $\alpha$ by translations parallel to the subtorus spanned by $x_4,\dots,x_7$ and due to the choice of neighborhoods $U$, their actions on $T^7/\langle \alpha \rangle$ are naturally extended to
actions on $T^3 \times K_3$.

The $\langle \beta,\gamma\rangle$-orbit of every fixed point torus of $\alpha$ consists of four tori.
Let us take the representatives  $T^3_{\alpha k}, k=1,2,3,4$, of these orbits.

Let us, without loss of generality, assume that $k=1$.

The domain $\widetilde{U_0/\langle \sigma \rangle}$ contains an embedded submanifold diffeomorphic to
$\cc\pp^1$ which realizes a nontrivial cycle $C_{\alpha 1}$. Let us take a closed two-form $c_{\alpha 1}$
 which is compactly supported in  $\widetilde{U_0/\langle \sigma \rangle}$ and realizes the Thom class of the normal bundle to $\cc\pp^1$ in $\widetilde{U_0/\langle \sigma \rangle}$. The pullback of this form onto
 $T^3_{\alpha 1} \times \widetilde{U_0/\langle \sigma \rangle}$ we define by the same symbol and
 extend it onto $T^3 \times K_3$ by zero.   It would be the Thom class of the normal bundle to $C^\prime_{\alpha 1}$. It is clear that the form $c_{\alpha 1}$ can be taken $\langle \beta,\gamma\rangle$-invariant.

To define the forms $c_{\alpha k}$ for $k=2,3,4$, i.e., for other orbits, we notice that all orbits are obtained from each other by translations by half-periods. We recall that by (\ref{comm})
all translations by half-periods commute with elements of $\Gamma$. Therefore, we obtain the forms
$c_{\alpha k}$ for $k=2,3,4$ by acting on $c_{\alpha 1}$ by convenient translations.

By construction, the forms $c_{\alpha k}$ are the pullbacks of the forms $s_{\alpha k}$ on $K_3$.
Since
$$
\int_{\widetilde{U_{k}/\langle \sigma \rangle}}c^2_{\alpha k} =
\int_{\widetilde{U_{j}/\langle \sigma \rangle}}c^2_{\alpha j}
$$ for all $j,k$ we have
$$
\int_{K_3} (s^2_{\alpha k} - s^2_{\alpha 1}) = \int_{t^\prime_\alpha} (c^2_{\alpha k} - c^2_{\alpha 1}) = 0
$$
and therefore the form $c^2_{\alpha k} - c^2_{\alpha 1}, k=2,3,4$,
is cohomologous to zero and there exists a form $\theta_{\alpha k}$ on $K_3$ such that
$$
d \theta_{\alpha k} = s^2_{\alpha k} - s^2_{\alpha 1} \ \ \mbox{on $K_3$}
$$
and
$$
d \eta_{\alpha k} = c^2_{\alpha k} - c^2_{\alpha 1} \ \ \ \mbox{on $\widetilde{T^7/\langle \alpha\rangle}$}
$$
for $k=2,3,4$.
In the domains where $c_{\alpha k}$ are supported, the forms $\eta_{\alpha j}$ and $c_{\alpha k}$
are linear combinations of the wedge products of the differentials $dx_4,dx_5,dx_6$, and $dx_7$
which correspond to  local coordinates on $\widetilde{U /\langle \alpha \rangle}$.

Let us consider the action of $\langle \beta, \gamma\rangle$ on $\widetilde{T^7/\langle \alpha \rangle}$.
The forms $c_{\alpha k}$ are invariant  with respect to this action, vanish near singular point sets of it and
induce the forms, on  $M$, which we denote also by $c_{\alpha k}$.
By analogy, we construct forms $c_{\beta j}$ and $c_{\gamma l}$ for $j,l=2,3,4$.

The averaged forms
$$
\eta^\prime_{\alpha k} = \frac{1}{4} \sum_{g \in \langle \beta,\gamma\rangle} g^\ast \eta_{\alpha k}
$$
also satisfy the equation
$$
d \eta^\prime_{\alpha k} = c^2_{\alpha k} - c^2_{\alpha 1} \ \ \ \mbox{on $\widetilde{T^7 /\langle \alpha \rangle}$}.
$$

Near the fixed point tori of $\beta$ and $\gamma$ the form $\eta^\prime_{\alpha k}$ is closed.
Moreover on a neihborhood $T^3 \times U$ of such a torus this form is cohomologous to zero.
Indeed, let, for instance, the torus be fixed by $\beta$, Then it corresponds to the coordinates
$x_3,x_6,x_7$ and the integral of $\eta^\prime_{\alpha k}$ over the torus vanishes. Let us take a smooth function on $U$ which is supported inside $U$, equal to $1$ on $U_0$, and invariant with respect to $\sigma= \beta$. Let us also take two-form $\omega$ in $T^3 \times U$ such that
$d\omega = \eta^\prime_{\alpha k}$, it is $\sigma$-invariant (which is easy to achieve by using the averaging).  Let us now replace $\eta^\prime_{\alpha k}$ in $T^3 \times U$ by the form
$$
\eta^{\prime\prime}_{\alpha k} = \eta^\prime_{\alpha k} - d(f \omega).
$$
Such a procedure can be done near every fixed torus of $\beta$ and $\gamma$ and, by the choice of $U_0$ and $U$, it can be done invariantly.
As a result we obtain a form $\eta^{\prime\prime}_{\alpha k}$ which is $\langle \beta, \gamma\rangle$-invariant and
vanishes near fixed point sets of the involutions $\beta,\gamma$, and $\beta\gamma$.
Therefore $\eta^{\prime\prime}_{\alpha k}$ induces a form on
$(\widetilde{T^7/\langle \alpha\rangle})/\langle \beta,\gamma\rangle$ which vanishes near the singular set and pushed back to the desingularization $M$ of this orbifold.
Therefore it gives rise two a form $n_{\alpha k}$ such that

1) $\int_M n_{\alpha k} c^2_{\delta j} = 0$ for $\delta  \neq \alpha$ because these forms supported on disjoint sets;

2) $\int_M n_{\alpha k} c^2_{\alpha j} = 0$ because on the subsets where these forms are supported
both forms are expressed by wedge products of $dx_4,dx_5,dx_6$, and $dx_7$ and these products vanish for dimension reasons.

Hence
\begin{align}
\label{eqn-added}
\int_M n_{\alpha i} c^2_{\delta k} = 0 \ \ \mbox{for all $i$ and $k$}.
\end{align}

By construction, we have
$$
dn_{\alpha k} = c^2_{\alpha  k} - c^2_{\alpha 1}.
$$

Therefore we construct the forms $c_{\delta i}$ and $n_{\delta k}$ with the desired prorepties.


\bigskip

Let us summarize in view of the previous constructions.
Degree $0$ is generated by the constant functions. In degree $3$, the subspace $\langle t_\delta, t_i, c_{\delta i j}\rangle_{\delta, i, j}$ just consists of representatives of the respective cohomology classes. The forms $n_{\delta k}$ and $n_\delta$ are chosen to satisfy
\begin{align}\label{eqn04}
\dif n_{\delta k}&= c_{\delta k}\wedge c_{\delta k}-c_{\delta 1}\wedge c_{\delta 1}\\
\label{eqn07}
\dif n_\delta&=c_{\delta i}\wedge c_{\delta i} -t_\delta'
\end{align}
for $2\leq k\leq 4$ representing the cohomological equalities $[c_{\delta i}]^2=[t_\delta']$ (and the thereby induced mutual equalities amongst the squares) for all $1\leq i\leq 4$. Note that by our choice of representatives in degree $2$ it holds that
\begin{align}\label{eqn03}
c_{\delta_1 i_1}\wedge c_{\delta_2 i_2}=0
\end{align}
for $\delta_1\neq \delta_2$ or $i_1\neq i_2$.

As a consequence, by construction, we obtain that the inclusion of cochain algebras $\iota\co (\mathcal{A},\dif) \to \ADR(M)$ is a quasi-isomorphism, i.e.~$H(\iota)$ yields the isomorphism
\begin{align*}
H(\mathcal{A},\dif)\xto{\cong} H(\ADR(M))\cong H^*(M)
\end{align*}

\bigskip

{\noindent \sc Passing to the minimal model.}

We now construct a minimal Sullivan model $(\Lambda V,\dif)$ (in normal form) for $(\mathcal{A},\dif)$ and thereby, by uniqueness of minimal models up to isomorphism, also doing so for $M$, i.e.~we construct a quasi-isomorphism
\begin{align*}
m\co (\Lambda V,\dif)\xto{\simeq} (\mathcal{A},\dif)\simeq \ADR(M)
\end{align*}
inductively over degree. That is, we construct the models and the morphisms
\begin{align*}
m^l\co (\Lambda V^{\leq l},\dif)\xto{} (\mathcal{A},\dif)\simeq \ADR(M)
\end{align*}
by constructing $V^{\leq l}$ together with the differential $\dif$ and the morphisms $m^l$ such that differentials commute with the morphism, and such that $m$ is a quasi-isomorphism up to degree $l$.

Set $V^1=0$, and obtain again by a little abuse of notation (and $m^2$, $m^3$ reflecting it exactly) that
\begin{align*}
V^2&=\langle c_{\delta i} \rangle_{\delta, i}\\
V^3&=\langle t_\delta, t_i, c_{\delta i j}\rangle_{\delta, i, j} \oplus \langle n_{\delta k}\rangle_{\delta, k} \oplus \langle q_s\rangle_s \\
\end{align*}
with the induced differentials as in \eqref{eqn04}. This procedure then may be continued over degree to yield $(\Lambda V,\dif)$ and $m$ in the limit as described in Algorithm \cite[Page 144]{FHT01}.
That is, we have constructed the quasi-isomorphisms
\begin{align*}
m\co (\Lambda V,\dif)\xto{\simeq} (\mathcal{A},\dif)\hto{\simeq} \ADR(M)\cong H^*(M)
\end{align*}

We need to comment on degree $3$ and, especially, on the additional summand $\langle q_s\rangle_s$: Since, in contrast to $(\mathcal{A},\dif)$, the algebra $(\Lambda V,\dif)$ is a free commutative differential graded algebra, it needs to encode the additional relations \eqref{eqn03} via differentials. Hence we choose the $q_s$ in bijection with these relations and set their differentials to equal each respective one. It is clear that $m$ can be set to restrict to an isomorphism on $\langle t_\delta, t_i, c_{\delta i j}\rangle_{\delta, i, j} \oplus \langle n_{\delta k}\rangle_{\delta, k}$ and can be set to vanish on $\langle q_s\rangle_s$.

Note that $n_\delta$ from Equation \eqref{eqn07} actually do not have to be represented in the minimal model, since the $t_\delta'$ just do not appear as new spherical cohomology; i.e.~in the minimal model---exactly due to minimality---they directly appear in the form $c_{\delta
 i}^2$ (for any $i$).

\bigskip

\subsection{3-formality of the model}
In view of Theorem \ref{theo01} it remains to show that this model $(\Lambda V,\dif)$ is $3$-formal. For this we split $V=C\oplus N$ as indicated above in such a way that
\begin{align*}
C^2&=\langle c_{\delta i} \rangle_{\delta, i} \\
N^2&=0\\
C^3&=\langle t_\delta, t_i, c_{\delta i j}\rangle_{\delta, i, j}\\
N^3&=\langle n_{\delta k }\rangle_{\delta, k} \oplus \langle q_s\rangle_s
\end{align*}

It remains to show that every closed form in
\begin{align}\label{eqn08}
\big(\langle n_{\delta k}\rangle_{\delta, k} \oplus \langle q_s\rangle_s\big)  \cdot (\Lambda V^{\leq 3})
\end{align}
is exact in $(\Lambda V,\dif)$. By degree and dimension this is only necessary in degrees $5$ to $7$, and by Poincar\'e duality this is clear in degree $6$. Hence we need to verify this in degrees $5$ and $7$.

\subsubsection{Degree $5$}
Suppose that there is a closed non-exact form from \eqref{eqn08} in degree $5$. By Poincar\'e duality, i.e.~by the non-degeneracy of the intersection form, there exist a cohomology class in degree $2$ which multiplies with the class in degree $5$ to a volume form. Consequently, this becomes a closed non-exact form from \eqref{eqn08} now in degree $7$. This implies that once we have proved that there are now such classes in degree $7$ the case of degree $5$ is settled as well.

We point the reader to the appendix where we have a closer look at dimension $5$, in particular illustrating the used techniques once again.

\bigskip

\subsubsection{Degree 7} \label{subsec7}
Now we show that there are no closed non-exact forms from \eqref{eqn08} in degree $7$ (thereby, as depicted, proving the analogue for degree $5$).

So we need to show that any closed form in $N^{\leq 3}\cdot (\Lambda V^{\leq 3})$ of degree $7$ is exact. We will do so by actually showing that \emph{any} form in there has an exact image in the quasi-isomorphic algebra $(\mathcal A,\dif)$.

By the structure of $V$ and the decomposition $V=C\oplus N$, it follows that any form in $N^{\leq 3}\cdot (\Lambda V^{\leq 3})$ of degree $7$ is necessarily in
\begin{align}
C^2\cdot C^2\cdot N^3
\end{align}

Note that by slight abuse of notation we will suppress the morphism $m$ in the following, etc., and consider forms to lie in the minimal model, the algebra $\mathcal{A}$ and the differential forms $\ADR(M)$ depending on context.

Next, in the algebra $(\mathcal A,\dif)$ it holds that $c_{\delta_1 i}\cdot c_{\delta_2 j}=0$ unless $\delta_1=\delta_2$ and $i=j$, since the forms are concentrated around disjoint submanifolds. Hence any non-exact element in $m(C^2\cdot C^2\cdot N^3)\In (\mathcal{A},\dif)$ (which is closed by degree) is necessarily in the vector space
\begin{align}\label{eqn06}
&\langle c_{\delta_1 i}^2\cdot n_{\delta_2  k}\rangle_{\delta_1, \delta_2, i, k}\In \mathcal{A}^7.
\end{align}

By \eqref{eqn-added} we now show that
\begin{align}\label{eqn02}
0=\int_{C'_{\delta_1 i}} \iota_{\delta_1 i}^*  (c_{\delta_1 i}\wedge n_{\delta_2 k})=\int_M (c_{\delta_1 i}\wedge n_{\delta_2 k})\wedge c_{\delta_1 i}
\end{align}
for the embeddings $\iota_{\delta_1 i}\co C'_{\delta_1 i} \to M$.
This proves that all the forms in \eqref{eqn06} above are exact. Hence $M$ is $3$-formal and formal consequently.

\begin{rem}
On the one hand in \cite{DGMS} it is stated that formality is equivalent to ``uniform vanishing of Massey products''. Since then a lot of work has been dedicated to making this statement precise in the form of several interpretations. This has led to various different characterizations of formality ranging from say the ``original'' one in \cite[Theorem 4.1, p.~261]{DGMS} over $A^\infty$- or $L^\infty$-algebras, etc.,~to the depiction cited in Remark \ref{rem03} below. On the other hand intersection theory was used before to represent Massey products via submanifolds. Hence in view of Theorem \ref{theo01} from \cite{FM05} which generalizes the characterization in \cite{DGMS} and which is our main algebraic tool we make this interplay of Massey products, formality, and submanifolds precise and geometric.
\end{rem}

\begin{rem}\label{rem04}
The reader may notice that our arguments, the blow-up construction, etc., do not really depend on the concrete structure of the group $\Gamma$. Hence this group certainly can be replaced by more general ones. Simply-connected manifolds are known to be formal below dimension $7$, whence also here a more general situation involving more complicated submanifolds should not lead to many additional difficulties. Hence, in summary, we are confident that our approach should pave the way to analyzing the formality of the remaining Joyce examples respectively of further instances like examples of $\Spin(7)$-manifolds, etc.
\end{rem}

\begin{rem}
We note that although the general theory (see for example \cite[Theorem 8.28, p.~337]{FOT08} and \cite{LS05}) for blow-ups allows for explicit rational models only in stable ranges (which in our case of a blow-up of a $3$-dimensional submanifold begins in dimension $9$) we managed to construct a rational model of $M$, which, due to formality clearly is given by its cohomology algebra. Hence, for example, it is a triviality to compute its rational homotopy groups (which clearly have exponential growth) in a certain range.
\end{rem}

\begin{rem}\label{rem03}
We remark that in the case of a simply-connected $7$-dimensional manifold, we may equivalently consider the ``Bianchi--Massey tensor'' introduced in \cite{CN}. Our arguments transcribe to show its vanishing, and according to \cite[Theorem 1.3, p.~3 ]{CN} this yields formality.
\end{rem}



\def\cprime{$'$}

\pagebreak \

\vfill

\begin{center}
\noindent
\begin{minipage}{\linewidth}
\small \noindent \textsc
{Manuel Amann} \\
\textsc{Institut f\"ur Mathematik}\\
\textsc{Differentialgeometrie}\\
\textsc{Universit\"at Augsburg}\\
\textsc{Universit\"atsstra\ss{}e 14 }\\
\textsc{86159 Augsburg}\\
\textsc{Germany}\\
[1ex]
\footnotesize
\textsf{manuel.amann@math.uni-augsburg.de}\\
\textsf{www.uni-augsburg.de/de/fakultaet/mntf/math/prof/diff/team/dr-habil-manuel-amann/}
\end{minipage}
\end{center}

\vspace{5mm}

\begin{center}
\noindent
\begin{minipage}{\linewidth}
\small \noindent \textsc
{Iskander A. Taimanov} \\
\textsc{Sobolev Institute of Mathematics}\\
\textsc{avenue academician Koptyug 4}\\
\textsc{630090 Novosibirsk}\\
\textsc{Russia}\\
[1ex]
\footnotesize
\textsf{taimanov@math.nsc.ru}
\end{minipage}
\end{center}

\end{document}